\documentclass[twoside,notitlepage,12pt]{article}

\usepackage{amssymb}
\usepackage[leqno]{amsmath}
\usepackage{amsfonts}
\usepackage{amsopn}
\usepackage{amstext}
\usepackage{amsthm}
\usepackage{enumerate}

\usepackage[colorlinks]{hyperref}
\usepackage{fourier-orns}
\usepackage{calrsfs}

\usepackage{makeidx}

\newcommand{\Ef}{{\cal{F}}}
\newcommand{\Gee}{{\cal{G}}}

\newcommand{\El}{{\cal{L}}}

\newcommand{\Tau}{{\cal{T}}}
\newcommand{\Yu}{{\cal{U}}}

\newcommand{\al}{\alpha}

\renewcommand{\phi}{\varphi}
\renewcommand{\rho}{\varrho}

\newcommand{\rest}{\restriction}

\newcommand{\ntr}{n\in\omega}

\newcommand{\loe}{\leqslant}

\newcommand{\subs}{\subseteq}
\newcommand{\sups}{\supseteq}
\newcommand{\nnempty}{\ne\emptyset}

\newcommand{\ovr}{\overline}
\renewcommand{\iff}{\Longleftrightarrow}
\newcommand{\iffdef}{\;\; \equiv \;\;}

\newcommand{\cl}{\operatorname{cl}}

\newcommand{\conv}{\operatorname{conv}}

\newcommand{\poset}{{\Bbb{P}}}

\newcommand{\meet}{\wedge}
\newcommand{\Meet}{\bigwedge}
\newcommand{\join}{\vee}

\newcommand{\Land}{\;\&\;}

\newcommand{\oraz}{\quad \text{ and } \quad}

\newtheorem{tw}{Theorem}[section]
\newtheorem{wn}[tw]{Corollary}
\newtheorem{lm}[tw]{Lemma}
\newtheorem{prop}[tw]{Proposition}

\theoremstyle{definition}

\newtheorem{ex}[tw]{Example}
\newtheorem{pyt}[tw]{Question}
\theoremstyle{remark}
\newtheorem{uwgi}[tw]{Remark}

\newcommand{\setof}[2]{\{#1\colon #2\}}

\newcommand{\seq}[1]{\langle #1 \rangle}
\newcommand{\seqof}[2]{\langle #1\colon #2\rangle}
\newcommand{\sett}[2]{\{#1\}_{#2}}
\newcommand{\sn}[1]{\{#1\}} % singleton
\newcommand{\dn}[2]{\{#1,#2\}} % doubleton
\newcommand{\pair}[2]{\langle #1, #2 \rangle} % pair
\newcommand{\triple}[3]{\langle #1, #2, #3 \rangle} % triple
\newcommand{\map}[3]{#1\colon #2 \to #3} % A function
\newcommand{\img}[2]{#1[#2]} % image of a set
\newcommand{\dpower}[2]{[#1]^{#2}}

\newcommand{\fin}[1]{[#1]^{<\omega}}

\providecommand{\cal}{\mathcal}
\renewcommand{\Bbb}{\mathbb}

\newenvironment{pf}{\begin{proof}}{\end{proof}}

\newcommand{\fs}{\operatorname{FS}}

\providecommand{\nat}{\omega}

\newcommand{\abs}[1]{|#1|}

\newcommand{\iso}{\cong}

\newcommand{\til}{\tilde}

\newcommand{\R}{\ensuremath{\mathcal R}}

\newcommand{\Szero}{\ensuremath{\mathbb S_0}}

\newcommand{\beto}{\mathfrak b_0}

\newcommand{\bb}{\mathbb}

\newcommand{\bD}{{\bb D}}
\newcommand{\bE}{{\bb E}}
\newcommand{\bG}{{\bb G}}
\newcommand{\bX}{{\bb X}}
\newcommand{\bY}{{\bb Y}}

\newcommand{\cL}{\ensuremath{\mathcal L}}
\newcommand{\cR}{\ensuremath{\mathcal R}}

\newcommand{\eva}[1]{\kappa_{#1}} % EVAluation map
\newcommand{\Eva}{\kappa} % EVAluation map

\newcommand{\X}{\ensuremath{\mathbb X}}
\newcommand{\Y}{\ensuremath{\mathbb Y}}
\newcommand{\Z}{\ensuremath{\mathbb Z}}
\newcommand{\K}{\ensuremath{\mathbb K}}
\newcommand{\bm}{\ensuremath{\mathbb M}}
\newcommand{\bd}{\ensuremath{\mathbb D}}
\newcommand{\bL}{\ensuremath{\mathbb L}}

\newcommand{\bO}{\ensuremath{\mathbb O}}

\newcommand{\wera}{\ensuremath{\mathbb W_{+0}}}
\newcommand{\vera}{\ensuremath{\mathbb W}}
\newcommand{\verazo}{\ensuremath{\mathbb W_{+01}}}
\newcommand{\werazo}{\verazo}

\newcommand{\bea}{\looparrowright}

\newcommand{\cmp}{\circ} % composition!!!
\newcommand{\univ}[1]{\ensuremath{\left|{#1}\right|}}
\newcommand{\isp}[1]{{\ensuremath{\operatorname{ISP}\!\left({#1}\right)}}}
\newcommand{\tisp}[1]{{\ensuremath{\operatorname{TISP}\!\left({#1}\right)}}}

\renewcommand{\hom}{\operatorname{Hom}}

\newcommand{\separator}{\begin{center}***\end{center}}

\title{%A survey of %KP
Dualities between finitely separated structures}
\author{
{\sc Wies{\l}aw Kubi\'s}\\ \\
Academy of Sciences of the Czech Republic\\
\texttt{kubis@math.cas.cz}
\and 
\mbox{}\\{\sc Krzysztof Pszczo{\l}a}\\ \\
Jan Kochanowski University\\
\texttt{pszczola@ujk.edu.pl}
}

\newcommand{\define}[1]{\index{#1}\emph{#1}}

\makeindex

\begin{document}
\maketitle

\begin{abstract} 
We study dualities between classes of relational topological structures, given by $\hom$-functors. 
We show that there exists a $2$-element structure with infinitely many relations, which reconstructs all other structures generated by a $2$-element one.
As an application, we find a natural duality for the class of normal convexity structures. 
As another application, we give short proofs for several known dualities for classes of structures generated by a fixed $2$-element structure.
%KP 

\ 

\noindent
{\bf MSC (2010)}
Primary:
03C52, %(1980-now) Properties of classes of models
08C20. %(2010-now) Natural dualities for classes of algebras
Secondary:
03E75, %(1980-now) Applications of set theory
06D50, %(2000-now) Lattices and duality 
18A40, %(1973-now) Adjoint functors (universal constructions, reflective subcategories, Kan extensions, etc.)
54H10. %(1973-now) Topological representations of algebraic systems

\noindent
{\bf Keywords and phrases:}
Semi-dual pair, duality, reflexivity, bi-convexity.
\end{abstract}

\tableofcontents

\section{Introduction}

By a \define{structure} we mean a set endowed with some constants, relations and operations (functions).
Many classes of structures appearing in general algebra and model theory are ``generated" by a fixed finite structure $\bD$ in the sense that, up to isomorphisms, the class consists of all substructures of arbitrary powers of $\bD$.
This is the case, for instance, with partially ordered sets, distributive lattices, semilattices, 2-groups, median algebras, etc.
We shall say that $\bX$ is \define{$\bD$-separated} if it embeds into some power of $\bD$.
This just means that the structure of $\bD$ is determined by homomorphisms into $\bD$.
Once we have a $\bD$-separated structure $\bX$, it is natural to look at $\hom(\bX,\bD)$, the set of all homomorphisms from $\bX$ to $\bD$.
Quite often, one can ``recognize" some (possibly different) structure on $\hom(\bX,\bD)$ which may happen to be $\bE$-separated for another finite structure $\bE$.
The best situation is when $\bE$ has the same universe as $\bD$ and the structure on $\hom(\bX,\bD)$ is inherited from the power $\bE^{|\bX|}$, where $|\bX|$ stands for the universe of $\bX$.
We write $\bX^* = \hom(\bX,\bD)$ endowed with the $\bE$-structure.
Now it is natural to consider $\bX^{**} = \hom(\bX^*, \bE)$ as a subset of $\bD^{\hom(\bX,\bD)}$.
There is a natural \define{evaluation map} $\map \Eva \bX {\bX^{**}}$ defined by $\Eva(x)(f) = f(x)$.
This map is always an embedding; if it is an isomorphism, we say that $\bX$ is \define{reflexive} 
(or -- more precisely -- reflexive with respect to $\bd$ and $\bE$). %KP
We talk about \define{natural duality} if all finite structures are reflexive.
We show that in this case all compact as well as all discrete structures are reflexive. 
(Note that the definition of \textit{natural duality} here is different than in \cite{CD} and other papers dealing with algebraic dualties. More on this can be found in Section \ref{SecNotationalremarks}.) %KP
In order to formulate the main statements, we need to talk about topological structures.

Namely, a structure is \define{topological} if the set has a Hausdorff topology such that all relations and all operations are continuous.
A structure $\bX$ is \define{$\bD$-separated} if it is embeddable into a power of $\bD$.
We now change the definition of $\hom(\bX,\bD)$ to the set of all \emph{continuous} homomorphisms from $\bX$ into $\bD$ (endowed with the discrete topology).
We may still ask for reflexive structures.
In the presence of natural duality, we show that the evaluation map $\map \Eva \bX {\bX^{**}}$ is surjective whenever $\bX$ is compact or discrete.

One of the best examples is the (nowadays classical) Priestley duality~\cite{DP}.
It says that bounded distributive lattices are reflexive with respect to the duality induced by the 2-element lattice and the 2-element partially ordered set.
Our approach gives ``for free" the reversed duality, discovered first by Gehrke (unpublished) and Banaschewski~\cite{Banasch}, much later reproved in \cite{ABKR}: partially ordered sets are dual to compact 0-dimensional distributive lattices.
It is interesting that a non-trivial topological ingredient is hidden here.
Namely, how to recognize that a certain compact 0-dimensional topological structure is topologically $\bE$-separated, where $\bE$ is the finite generating structure?
We say that a compact structure is \define{topologically $\bE$-separated} if continuous homomorphisms into $\bE$ separate points.
It is proved in \cite{HMS} (and also reproved in \cite{ABKR}) that this happens with \emph{all} distributive lattices (here $\bE$ is the 2-element lattice).
It is easy to see that this does not happen with compact partially ordered spaces, the $2$-separated ones are called \define{Priestley spaces}~\cite{DP}.

Extensive study of natural dualities is given in the monograph of Clark \& Davey~\cite{CD}, where the authors consider algebras and allow relations in compact structures only.
These assumptions have been relaxed in several recent papers \cite{Hof}, \cite{Davey2007}, \cite{DHP}; for applications see also \cite{Jo}. %KP

We propose a unified approach, emphasizing on relational structures.
In particular, we describe an ultimate structure with the two-element universe that encodes all 2-separated structures and offers a useful duality whose special cases are several well known dualities.
We also discuss some applications to abstract convexity structures.

%abcd

\section{Natural dualities revisited}\label{SecNatDual}

Let $\bm=\seq{M,\sett{R_i}{i\in I},\sett{f_j}{j\in J},\sett{c_k}{k\in K}}$ be a model of a fixed first-order language $\El$, where $R_i$, $f_j$ and $c_k$ are predicates for relations, (partial) operations and constants, respectively.
Let $\bb X$ be another structure of the same language $\El$. We say that $\bb X$ is \define{$\bm$-separated} if there is an embedding of $\bb X$ into some power of $\bm$.
By an \define{embedding} we mean an injective homomorphism $\map f\X\Y$ satisfying the equivalence
$$R(f(x_0),\dots,f(x_{n-1})) \iff R(x_0,\dots,x_{n-1})$$
for every $n$-ary relation $R$ in the common language of $\X$ and $\Y$.
In other words, $f$ is an embedding iff the image $\img f\X$ is a substructure of $\Y$ isomorphic to $\X$ via $f$.

Following~\cite{CD}, we denote by $\isp{\bm}$ the class of all $\bm$-separated structures.
It can be described as the class of all structures isomorphic to substructures of arbitrary powers of $\bm$ (hence the shortcut ISP).
In fact, given an $\bm$-separated structure $\X$, there is an embedding $\map j \X {\bm^H}$, given by $j(x)(h) = h(x)$, where $H = \hom(\X,\bm)$.
On the other hand, every substructure of $\bm^\kappa$ is clearly $\bm$-separated.

In case where $\bm$ is finite, it can be regarded as a discrete topological structure and its powers are compact 0-dimensional structures.
We denote by $\tisp \bm$ the class of all topological structures in the language of $\bm$, admitting a continuous embedding into some power of $\bm$.
When dealing with topological structures, the set $\hom(\bX,\bY)$ will always denote all continuous homomorphisms from $\bX$ to $\bY$.

Given a model $\bm$, we denote its universe by $\univ\bm$. We say that $\bm$ is {\em finite} if $\univ\bm$ is a finite set.
We shall write $X,Y,Z, \dots$ for the universe of $\X,\Y,\Z,\dots$.

\subsection{Dual pairs}

Fix two countable first-order languages $\cL$, $\cR$ and fix two finite models $\bD$, $\bE$ of $\cL$, $\cR$ respectively, so that $\univ \bD = \univ \bE$.
We say that $\pair \bD\bE$ is a \define{semi-dual pair} if for every finite $\X\in \isp\bD$  
the following conditions are satisfied:
%KP
\begin{enumerate}
	\item[(S1)] The set $\hom(\X,\bD)$ is a substructure of the power $\bE^X$.
	\item[(S2)] For every $\phi\in\hom(\hom(\X,\bD),\bE)$ there exists $x\in X$ such that $\phi(f)=f(x)$ for $f\in \hom(\X,\bD)$. Here we consider the set $\hom(\X,\bD)$ with the $\cR$-structure induced from $\bE^X$.
\end{enumerate}
Finally, $\pair \bD\bE$ is a \define{dual pair} if both $\pair\bD\bE$ and $\pair \bE\bD$ are semi-dual.

Assume $\pair\bD\bE$ is a fixed dual pair.
We shall use the following abbreviations.
Given a topological model $\X\in\tisp\bD$, we shall denote by $\X^*$ the set $\hom(\X,\bD)$ endowed with the $\cR$-structure induced from the power $\bE^X$ (recall that we consider continuous homomorphisms only).
Notice that $\hom(\X,\bD)$ is indeed a substructure of $\bE^{|\X|}$, because of condition (S1).
Indeed, given a $k$-ary function symbol $F$ in $\cR$, given $f_1,\dots,f_k \in \hom(\X,\bD)$, the map $F(f_1,\dots,f_k)$ is defined pointwise, therefore it is a homomorphism (because its restriction to every finite substructure of $\X$ is a homomorphism by (S1)) and it is continuous with respect to the product topology.

Similarly, given a topological structure $\Y\in \tisp\bE$, we denote by $\Y^*$ the $\cL$-structure on $\hom(\Y,\bE)$ induced from $\bD^Y$, where $Y=\univ\Y$.
Finally, given $\X\in \tisp\bD$, we identify $x\in X$ with the homomorphism $\eva x\in\hom(\X^*,\bD)$ defined by $\eva x(f)=f(x)$. Notice that the map $x\mapsto \eva x$, called the \define{evaluation map}, is indeed an embedding of $\X$ into $\X^{**}$.
It is moreover continuous with respect to the product topology on $\X^{**}$ which is just the topology of pointwise convergence.
The evaluation map is often not a topological embedding, simply because the topology of $\X$ may be strictly bigger than the weak topology induced by all continuous homomorphisms into $\bD$.

We call $\X^*$ and $\X^{**}$ respectively the \define{dual structure} and the \define{second dual structure} of $\X$ (relatively to $\pair\bD\bE$). If $\pair \bD\bE$ is a dual pair then, by definition, for every finite $\X\in\isp\bD$ the evaluation map $\map \Eva \X{\X^{**}}$ is onto, i.e. it is an isomorphism.

Now suppose that $\pair \bD\bE$ is a semi-dual pair only.
Given $\bX \in \tisp \bD$, define $\bX^*$ and $\bX^{**}$ as before, ignoring the fact that $\bX^{**}$ may not be a substructure of $\bD^{|\bX^*|}$ (to be formal, we may let $\bX^{**}$ to be generated by $\hom(\bX^*,\bE) \subs \bD^{|\bX^*|}$).
The evaluation map $\map \Eva \bX {\bX^{**}}$ is still an embedding and we can ask when it is an isomorphism.
If this is so, automatically the set $\hom(\bX^*,\bE)$ is a substructure of $\bD^{|\bX^*|}$.

Here is the formal definition.
A structure $\X$ will be called \emph{reflexive} with respect to a semi-dual pair $\pair \bD\bE$ if the evaluation map $\map\Eva \X{\X^{**}}$ is onto.
Condition (S2) in the definition of a semi-dual pair says that every finite structure is reflexive.
We are going to prove below that both compact and discrete structures are reflexive with respect to a fixed semi-dual pair.
We start with the first (easier) part.

\begin{tw}\label{tdwajedenxx}
Let $\pair \bD\bE$ be a semi-dual pair and let $\K \in \tisp \bD$ be a compact structure.
Then
\begin{enumerate}
	\item[(a)] Every homomorphism from $\K^*$ into $\bE$ is continuous.
	\item[(b)] $\K$ is reflexive 
with respect to $\pair \bD\bE$. %KP
\end{enumerate}
\end{tw}

\begin{pf}
We shall prove (a) and (b) simultaneously.
Let $\Y = \K^*$.
We may assume that $\K \subs \bD^\Y$.
Fix a finite set $S\subs \Y$ and define
$$\K(S) = \setof{ x\rest S}{x \in \K}.$$
Then $\K(S)$ is a finite $\bD$-separated structure.
Given $f\in \K(S)^*$, let $\til f (x) = f(x\rest S)$.
This defines an embedding of $\K(S)^*$ into $\Y$.
Let
$$\Y(S) = \setof{\til f}{f\in \K(S)^*}.$$
Then $\Y(S)$ is canonically isomorphic to $K(S)^*$.

Now fix a homomorphism $\map \phi {\K^*}\bE$.
We do not assume in advance that $\phi$ is continuous.
Again fix a finite $S\subs \Y = \K^*$.
Then $\phi\rest \Y(S) \in \Y(S)^*$ and $\Y(S) \iso \K(S)^*$, therefore
there exists $x_S\in \K(S)$ such that
$$\phi(f) = f(x_S)$$
for every $f\in \Y(S)$.
In particular, $x_S(s) = \phi(p_s)$, where $p_s\in \K(S)^*$ is the projection on the $s$-th coordinate, i.e. $p_s(x) = x(s)$ for $x\in \K$.
Thus, if $S\subs T\subs \Y$ then $x_S = x_T \rest S$.

Let $x\in \bD^\Y$ be such that $x\rest S = x_S$ for every finite $S\subs \Y$.
We claim that $x\in \K$.
In fact, given a finite $S\subs \Y$, there exists $y_S\in \K$ such that $y_S\rest S = x_S$.
This shows that $x$ belongs to the closure of $\K$.
It follows that $x\in\K$, because $\K$ is compact and hence closed in $\bD^\Y$.
In particular, we have shown that $\phi$ is continuous.
\end{pf}

The proof that discrete structures are reflexive requires more work.
We start with some lemmas.
The first one is quite well-known, we present a proof just for completeness.

\begin{lm}\label{lgenfinstru}
Let $\bD$ be a finite structure and let $\X$ be $\bD$-separated. Then every finitely generated substructure of $\X$ is finite.
\end{lm}

\begin{pf}
We may assume that $\X$ is a substructure of $\bD^\kappa$ for some (infinite) set $\kappa$.
Fix a finite $S\subs \X$ and let $S = \{ s_0,\dots, s_{n-1} \}$.
Define $\zeta(\al) = \seq{ s_0(\al), \dots, s_{n-1}(\al) } \in \bD^n$.
Note that $F := \img \zeta S \subs \bD^n$ is finite.
Write $\kappa = \bigcup_{\xi\in F}\kappa_\xi$, where $\kappa_\xi = \zeta^{-1}(\xi)$.
Notice that each element of $S$ is constant on every $\kappa_\xi$.
Now let $G$ consist of all $x\in \bD^\kappa$ that are constant on each $\kappa_\xi$.
It is clear that $G$ is finite; in fact $|G| = |\bD|^{|F|}$.
It is also clear that $G$ is a substructure of $\bD^\kappa$ and $S\subs G$.
\end{pf}

\begin{lm}\label{lewrjgpwe}
Let $\pair\bD\bE$ be a semi-dual pair and let $\map f\X\Y$ be a surjective homomorphism of $\bD$-separated structures. Then the dual map $\map {f^*}{\Y^*}{\X^*}$ defined by
$f^*(y^*) = y^* \cmp f$ for $y^*\in\Y^*$, is an embedding.
\end{lm}

\begin{pf}
Since $f$ is onto, it is clear that $f^*$ is one-to-one.
It is also clear that $f^*$ is a homomorphism of $\bE$-separated structures.
Fix a relation $R$ in the language of $\bE$ and suppose that $\neg R(y_1^*,\dots, y_k^*)$.
Then there exists $y\in \Y$ such that $\neg R(y_1^*(y), \dots, y_k^*(y))$.
Find $x\in \X$ so that $y = f(x)$.
Then $\neg R(f^*(y_1^*)(x), \dots, f^*(y_k^*)(x))$.
It follows that $\neg R(f^*(y_1^*), \dots, f^*(y_k^*))$.
This shows that $f^*$ is an embedding.
\end{pf}

\begin{lm}\label{lwtoiotu}
Assume $\pair \bD\bE$ is a semi-dual pair and let $\X$ be a discrete $\bD$-separated structure.
Then for every finite substructure $S$ of $\X$ there exists a finite structure $G(S)$ such that
\begin{enumerate}
	\item[(1)] $S\subs G(S) \subs \X$,
	\item[(2)] $G(S)^*$ is embeddable into $S^*$,
	\item[(3)] every homomorphism $\map f{G(S)}\bD$ extends to a homomorphism $\map {\ovr f}\X\bD$.
\end{enumerate}
\end{lm}

\begin{pf}
Fix a finite structure $Y\subs\X$ such that $S\subs Y$. Let $\map {r_Y}{Y^*}{S^*}$ be the restriction map and let $H(Y) = \img {r_Y}{Y^*}$.
We have two homomorphisms $\map r{Y^*}{H(Y)}$ and $\map i{H(Y)}{S^*}$, where the first one is $r_Y$ with restricted co-domain and the latter one is just the inclusion. Furthermore, $r_Y = i \cmp r$.
Using duality, we have that $j = r^* \cmp i^*$, where $j$ is the inclusion $S\subs Y$.
It follows that $\map{i^*}{S}{H(Y)^*}$ is an embedding.
Thus, we may assume that $S \subs H(Y)^* \subs Y$.
Clearly, every $f\in \hom(H(Y)^*,\bD)$ extends to some $\ovr f\in \hom(Y,\bD)$.

Now, if $Y_1\subs \X$ is such that $Y\subs Y_1$ then $H(Y_1) \subs H(Y)\subs S^*$.
Since $S^*$ is finite, there exists a finite $Z\subs \X$ such that $H(Y) = H(Z)$ whenever $Y$ is finite and such that $Z \subs Y \subs \X$.
We claim that $G(S) := H(Z)^*$ is as required.

Properties (1) and (2) are obvious.
Given a finite structure $Y\sups Z$, by the above arguments we know that every homomorphism $\map f{G(S)}{\bD}$ extends to a homomorphism $\map{f_Y}{Y}{\bD}$.
Let $\Ef$ be the family of all finite structures of the form $G(Y)$ where $Y\subs \X$ is finite and $Y\sups Z$.
Given $F_0\subs F_1$ in $\Ef$, we know that every $f_0\in F_0^*$ extends to some $f_1\in F_1$.

Fix $f\in G(S)^*$. For $F\in\Ef$ define
$$H_F = \setof{x\in \bD^\X}{x\rest F \in F^* \text{ and } x\rest G(S) = f}.$$
Notice that $H_F$ is a closed (in fact: clopen) subset of $\bD^\X$.
It is nonempty by the above arguments.
Given $F_0,\dots, F_{k-1} \in \Ef$ we have that
$$H_F \subs \bigcap_{i < k}H_{F_i},$$
where $F\in\Ef$ is such that $F_i\subs F$ for every $i < k$.
By compactness, there exists $y\in \bigcap_{F\in\Ef}H_F$.
Then $y\rest G(S) = f$ and clearly $y\in \X^*$.
\end{pf}

\begin{tw}\label{tdwaipull}
Let $\pair \bD\bE$ be a semi-dual pair and let $\X\in \isp\bD$ be a discrete structure.
Then $\X$ is reflexive with respect to $\pair \bD\bE$.
\end{tw}

\begin{pf}
Let $\K = \X^*$. Note that $\K$ is compact.
Given $S\subs \X$, define $\K(S) = \setof{y\rest S}{y\in \K}$.
Fix $\phi \in \K^*$.
Since $\bE$ is finite, using the compactness of $\K$, we can find a finite $G\subs \X$ such that $\phi$ depends on coordinates from $G$, i.e.
$$(\forall\; f_0,f_1\in \K)\; f_0 \rest G  = f_1 \rest G \implies \phi(f_0) = \phi(f_1).$$
By Lemma~\ref{lwtoiotu}, we may assume that $G = G(S)$ for some $S$, that is, the restriction map $\map r{\X^*}{G^*}$ is onto.
In other words, $G^* = \K(G)$.
Now observe that $\phi$ induces a unique homomorphism $\map{\til\phi}{G^*}{\bE}$ satisfying $\phi = \til\phi \cmp r$ or, in other words, $\til\phi(f\rest G) = \phi(f)$ for $f\in \K$.

Using reflexivity we can find $x\in G$ such that $\til\phi(g) = g(x)$ for every $g\in G^*$.
Finally, we get
$$\phi(f) = \til\phi(f\rest G) = f(x)$$
for every $f\in \K$.
This shows that $\phi$ is represented by a point of $\X$.
\end{pf}

%\fourierseparator

%The following statement might be false \grimace\grimace\grimace, although we do not have any candidate for a counterexample. (It seems that it is true when $\bD$ (or $\bE$) is injective.)

%\begin{tw}
%Let $\pair \bD\bE$ be a dual pair and let $\X\in \tisp\bD$ be such that $\X^*$ is compact with respect to the pointwise topology. Then $\X$ is reflexive with respect to $\pair \bD\bE$. Furthermore, every homomorphism from $\X$ into $\bD$ is continuous.
%\end{tw}

%\begin{pf}
%Let $\K = \X^*$ and fix a continuous homomorphism $\map \phi\K\bE$. Since $\bE$ is finite, there exists a finite set $S\subs \X$ such that $$(\forall\; f_0,f_1\in \K)\; f_0 \rest S  = f_1 \rest S \implies \phi(f_0) = \phi(f_1).$$ Enlarging $S$ if necessary, we may assume that it is a substructure of $\X$ (see Lemma~\ref{lgenfinstru}).??????
%\end{pf}

\begin{uwgi}
The main point in the results above is the full symmetry in case of dual pairs. Namely, given a dual pair $\pair \bD \bE$, we have in fact $4 = 2 \times 2$ theorems on the reflexivity of compact / discrete / $\bD$-separated / $\bE$-separated structures.
%Looking at the arguments above, it is easy to notice that actually we have only used the property that $\pair \bD \bE$ is a semi-dual pair. In fact, there is no need to formally define $\bX^{**}$, reflexivity just means that every element of $\hom(\bX^*,\bE)$ is induced by an element of $\abs \bX$. In other words, we do not need to know in advance whether the pair $\pair \bE \bD$ satisfies (S1). As the conclusion, we see that Theorems~\ref{tdwajedenxx} and~\ref{tdwaipull} are true under the assumption that $\pair \bD \bE$ is a semi-dual pair.
\end{uwgi}

%KP
%Here we summarize the most obvious implications of the results presented in this section:

%If $\pair \bD\bE$ is a dual pair, then we have the following dualities:
%\begin{itemize}
%\item Between the classes $\isp\bD_{\text{fin}}$ and $\isp\bE_{\text{fin}}$ (by definition of a dual pair; by $\mathcal{A}_{\text{fin}}$ we mean the subclass of $\mathcal{A}$ consisting of all finite elements of $\mathcal{A}$).
%\item Between the classes $\isp\bD$ and $\tisp\bE$,  and between the classes $\isp\bE$ and $\tisp\bD$ (by Theorem \ref{tdwajedenxx}).
%\item Between the classes $\isp\bD$ and $\isp\bE$ (by Theorem \ref{tdwaipull}).
%\end{itemize}
%KP

\begin{uwgi}
It is possible to prove Theorems~\ref{tdwajedenxx} and \ref{tdwaipull} using category-theoretic methods, more precisely, direct and inverse systems and their (co-)limits.
In fact, given a finite structure $\bD$, the category of finite $\bD$-separated structures is co-dense in the category of all discrete $\bD$-separated structures and dense in the category of all compact $\bD$-separated structures.
For definitions and some basic results concerning (co-)density of categories we refer to Chapter 0 of \cite{HMS} (see also \cite{Isbell} for a more general approach).

We have decided to present elementary arguments of rather set-theoretic nature, in order to make the results accessible to a more general audience.
Of course, by this way we are hiding some important ideas from category theory.

Just to summarize, given a dual pair $\pair \bD \bE$, the main results of this section say that the category of discrete $\bD$-separated structures is dually equivalent to the category of compact $\bE$-separated structures, and vice-versa.

Category-theoretic approach to natural dualities can be found in \cite{Hof}.
\end{uwgi}

\section{Dualities for 2-separated structures} \label{Sec2sep}

In this section we study dualities for $\bD$-separated structures, where $\bD$ is a structure whose universe is $2=\dn01$.
This includes posets, distributive lattices, semilattices and so on.
We show below that there exists a universal structure on $\dn 01$ which provides semi-duality with all other structures on $\dn 01$.

As it happens, the constants play a significant role for dualities.
Thus, we make the following agreement.
Given a structure $\bD$ whose universe is $\dn 01$, we shall denote by $\bD_{+0}$, $\bD_{+1}$ and $\bD_{+01}$ the same structure with constants $0$, $1$ or both of them, respectively, added to the language.

\subsection{An ultimate structure with the 2-element universe}\label{SubSecUlt}

Define $\vera = \pair{2}{\sett{I_{n,m}}{n,m \in\nat}}$, where $I_{n,m}$ is a relation of arity $n+m$ defined by
$$I_{n,m}(a_0, \dots, a_{n-1}, b_0, \dots, b_{m-1}) \iff \min_{i < n} a_i \loe \max_{j < m} b_j.$$
We allow $n=0$ or $m=0$, agreeing that $\min \emptyset = 1$ and $\max \emptyset = 0$.
In particular, $I_{0,0}$ is the empty relation and plays almost no role.
The following easy characterization partially explains the meaning of the relations $I_{n,m}$:

\begin{prop}\label{PropTrzideset}
Every $\vera$-separated structure is isomorphic to
$$\seq{\Ef, \seqof{\i_{n,m}}{n,m < \nat}},$$
where $\Ef$ is a family of sets, and
$$\i_{n,m}(A_0,\dots,A_{n-1},B_0, \dots, B_{m-1}) \iff \bigcap_{i < n}A_i \subs \bigcup_{j < m}B_j.$$
\end{prop}

\begin{pf}
Let $\X$ be a $\vera$-separated structure. We may assume that it is a substructure of $(\vera)^Y$, where $Y = \X^*$.
Identify $(\vera)^Y$ with the powerset of $Y$. Then $\X$ becomes a family of sets, the zero becomes the empty set.
Finally, $I_{n,m}(a_0,\dots,a_{n-1},b_0,\dots,b_{m-1})$ holds if and only if every homomorphism $\map h\X\vera$ such that $h(a_i)=1$ for every $i < n$, satisfies $h(b_j)=1$ for some $j < m$.
The last condition is equivalent to the one stated above, i.e. $I_{n,m} = \i_{n,m}$.
\end{pf}

From now on, we shall introduce a new (slightly non-standard) notation, namely, instead of the relations $I_{n,m}$ we shall work with a single relation $\bea$ defined on pairs of finite sets, having in mind that it encodes all the relations $I_{n,m}$.
More precisely, define
$$a \bea b \iff I_{n,m}(s_0, \dots, s_{n-1},t_0,\dots, t_{m-1})$$
for every $a = \{s_0, \dots, s_{n-1}\}$, $b = \{t_0, \dots, t_{m-1}\}$.
This new relation $\bea$ ``hides" the properties saying that the relations $I_{n,m}$ are stable under suitable permutations of parameters.
We can now list the set of axioms relevant for our structures.
\begin{enumerate}
\item[($\i_0$)] $\emptyset \not\bea \emptyset$.
\item[($\i_1$)] $a \bea b \Land a \subs a' \Land b \subs b' \implies a' \bea b'$.
\item[($\i_2$)] $\sn p \bea \sn q \Land \sn q \bea \sn p \iff p = q$.
\item[($\i_3$)] $(a_0 \cup \sn p) \bea b_0 \Land a_1 \bea (b_1 \cup \sn p) \implies (a_0 \cup a_1) \bea (b_0 \cup b_1)$.
\end{enumerate}
Condition ($\i_3$) will be called the \define{Pasch axiom}, its meaning will be explained later.
It is rather clear that all $\vera$-separated structures satisfy ($\i_0$)--($\i_3$).
Actually, the Pasch axiom requires some easy computations.
The converse is less obvious, as we see below.

\begin{lm}\label{LemBeaaeb}
Let $\bX$ be a $\vera$-separated structure and let $\dn L U$ be a partition of $X$.
Then the characteristic function of $U$ is a homomorphism if and only if
\begin{equation}
(\forall\; s \in \fin U)(\forall\; t \in \fin L) \;\; s \not \bea t.
\tag{\L}\label{Eqnotbaoab}
\end{equation}
\end{lm}

\begin{pf}
Let $f = \chi_U$.
If $f$ preserves $\bea$ then it necessarily satisfies (\ref{Eqnotbaoab}).
Suppose $f$ is not a homomorphism, that is, there are finite sets $s$, $t$ such that $t \bea s$ and $\img f t \not \bea \img f s$.
It follows that $\img f t \subs \sn 1$ and $\img f s \subs \sn 0$.
In other words, $s \subs L$ and $t \subs U$.
Thus, $f$ fails condition (\ref{Eqnotbaoab}).
\end{pf}

\begin{tw}\label{ThmPaschNormal}
Let $\bX = \pair X \bea$ be a structure, where $\bea$ is a binary relation on pairs of finite sets satisfying $(\i_0)$--$(\i_3)$.
Then $\bX \in \isp \vera$.
\end{tw}

\begin{pf}
Note that condition ($\i_2$) says that equality is defined from $\bea$.
Thus, in order to show that $\bX$ is $\vera$-separated, we need to show that, given finite sets $a, b \subs X$ satisfying $a \not \bea b$, there exists a homomorphism $\map h \bX \vera$ such that $\min \img h a = 1$ and $\max \img h b = 0$.
By ($\i_2$), this will also show that homomorphisms into $\vera$ separate points.

Fix $a,b$ with $a \not\bea b$ and let $U \sups a$ be a maximal (with respect to inclusion) set satisfying
\begin{enumerate}
	\item[(1)] $U \cap b = \emptyset$.
	\item[(2)] $(\forall\; s \in \fin U)\;\;s \not \bea b$.
\end{enumerate}
Clearly, such a set exists by Zorn's Lemma, using the fact that $a$ satisfies (1), (2) in place of $U$.
Now, let $L \sups b$ be a maximal set satisfying (\ref{Eqnotbaoab}) of Lemma~\ref{LemBeaaeb} and
\begin{enumerate}
	\item[(3)] $L \cap U = \emptyset$.
\end{enumerate}
Again, the existence of $L$ follows from Zorn's Lemma (note that $L := b$ satisfies (\ref{Eqnotbaoab}) and (3)).
We claim that $L \cup U = X$.

Suppose otherwise and fix $p \in X \setminus (L \cup U)$.
By the maximality $U$, the set $U \cup \sn p$ does not satisfy (2), so there is $s_0 \in U$ such that $(s_0 \cup \sn p) \bea b$.
By the maximality of $L$, the set $L \cup \sn p$ fails (\ref{Eqnotbaoab}), therefore there are $s_1 \in \fin U$, $t_1 \in \fin L$ such that $s_1 \bea (t_1 \cup \sn p)$.
The Pasch axiom (condition ($\i_3$)) tells us that $(s_0 \cup s_1) \bea (b \cup t_1)$, which contradicts (\ref{Eqnotbaoab}), because $b \subs L$.
Thus we have proved that $X = L \cup U$.

Finally, by Lemma~\ref{LemBeaaeb}, the characteristic function of $U$ is a homomorphism such that $\img h a \subs \sn 1$ from $\img h b \subs \sn 0$.
\end{pf}

It is rather easy to guess the additional axiom for the class of $\wera$-separated structures.
Namely, in view of Proposition~\ref{PropTrzideset}, the relevant axiom is
\begin{equation}
\sn 0 \bea \emptyset.
\tag{$c_0$}\label{EqC0}
\end{equation}
Similarly, for $\werazo$-separated structures we need to add (\ref{EqC0}) and
\begin{equation}
\emptyset \bea \sn 1
\tag{$c_1$}\label{EqC1}
\end{equation}
to the list ($\i_0$)--($\i_3$).
It is interesting that axioms (\ref{EqC0}), (\ref{EqC1}) together with ($\i_0$)--($\i_3$) allow us to define the notion of a complement (or negation) of an element.
We explain it later, when discussing some other features of the relation $\bea$.

We now show that $\werazo$ can be regarded as the ultimate structure on the two-element universe.

\begin{tw}\label{ThmDawidek}
Let $\bD$ be a structure whose universe is $\dn 01$.
Then $\pair \bD \werazo$ is a semi-dual pair.

Furthermore, if $0$ belongs to the language of $\bD$ then $\pair \bD \wera$ is a semi-dual pair.
Similarly, if both $0$ and $1$ are in the language of $\bD$ then $\pair \bD \vera$ is a semi-dual pair.
\end{tw}

\begin{pf}
Fix a finite $\bX \in \isp \bD$, where we assume that $|\bD| = \dn01$.
Fix $\phi \in \bX^{**}$.
Since $0,1$ are in the language of $\werazo$, $\phi$ is not constant.
Consider the universe of $\bX^*$ as the family of all subsets of $X$, whose characteristic functions are the homomorphisms into $\bD$.
Let $\Ef = \phi^{-1}(1)$, treated as a family of sets.
Note that $\Ef \nnempty$, because $\phi \ne0$.

Suppose $\bigcap \Ef = \emptyset$.
Then $\Ef \bea \emptyset$ and $\img \phi \Ef = \sn 1$, therefore $\phi$ is not a homomorphism, a contradiction.

Thus $\Ef \nnempty$ and $G := \bigcap \Ef \nnempty$.
Suppose for each $x \in G$ there is $f_x \in \bX^*$ such that $\phi(f_x) \ne f_x(x)$.
Notice that $F_x := f_x^{-1}(1) \notin \Ef$, because otherwise we would have that $f_x(x) = 1$ and also $\phi(f_x) = 1$.
Thus $\phi(f_x) = 0$, therefore $f_x(x) = 1$, that is, $x \in F_x$.
Finally, we have
$$\bigcap \Ef = G \subs \bigcup_{x \in G}F_x$$
which means that $\Ef \bea \sett{F_x}{x \in G}$.
This again shows that $\phi$ is not a homomorphism, a contradiction.

Finally, if $0$ is in the language of $\bD$, then we allow the possibility that $\phi = 0$ and in this case $\phi(f) = f(0)$ for $f \in \bX^*$.
If both $0$ and $1$ are in the language of $\bD$ then we also allow that $\phi = 1$ and consequently $\phi(f) = f(1)$ for $f \in \bX^*$.
This completes the proof.
\end{pf}

As a corollary, we obtain the main result of this section:

\begin{tw}\label{TwCorWeronisia}
$\pair \wera \wera$ and $\pair \vera \werazo$ are dual pairs.
\end{tw}

Given a structure $\bD$ whose universe is $\dn 01$, given a $\bD$-separated structure $\bX$, it is natural to define
\begin{equation}
s \bea t \iffdef (\forall\; h \in \hom(\bX,\bD)) \;\; \img h s \bea \img h t,
\tag{H}\label{Eqwationae}
\end{equation}
where $s, t \in \fin X$.
It is easy to see that such defined relation $\bea$ satisfies ($\i_0$)--($\i_3$), therefore $\bX$ becomes a $\vera$-separated structure.
It is natural to ask whether $\bea$-homo\-mor\-phisms into $\vera$ are the same as homomorphisms into $\bD$.
This turns out to be true, as we show below.

\begin{tw}\label{Twthmeorgo}
Let $\bD$ be a structure with universe $\dn 01$ and let $\bX \in \isp \bD$.
Consider $\bX$ with the relation $\bea$ defined by condition (\ref{Eqwationae}) above.
Then $$\hom(\bX,\bD) = \hom(\bX, \vera),$$
when treated just as sets with no structure.
\end{tw}

\begin{pf}
Just by definition, we have that $\hom(\bX,\bD) \subs \hom(\bX, \vera)$.
Fix $\map f X {\dn 01}$ which preserves $\bea$ and suppose that $f$ is not a homomorphism into $\bD$.
Replace each (possibly partial) operation in the language of $\bD$ by a suitable relation.
This causes no problem, because we deal with a fixed $\bD$-separated structure.
Now the fact that $f$ is not a homomorphism is witnessed by a relation $R$ in the language of $\bD$.
More precisely, there exist $a_0, \dots, a_{n-1} \in X$ such that
$$\bX \models R(a_0, \dots, a_{n-1}) \quad \text{ while }\quad \bD \not \models R(f(a_0), \dots, f(a_{n-1})).$$
Changing the order if necessary, we may assume that $f(a_i) = 0$ for $i < k$ and $f(a_i) = 1$ for $k \loe i < n$ (possibly $k=0$ or $k=n$).
Let $t = \setof{a_i}{i < k}$ and $s = \setof{a_i}{k \loe i < n}$.
Then $s \not \bea t$ in $\bX$, because $f$ preserves $\bea$ and $\img f s \not \bea \img f t$ in $\vera$.
By the definition of $\bea$ in $\bX$ (formula (\ref{Eqwationae})), there exists $h \in \hom(\bX, \bD)$ such that $\img h s \not \bea \img h t$.
This means that $\img h s = \img f s$ and $\img h t = \img f t$, therefore
$$\bD \not \models R(h(a_0), \dots, h(a_{n-1})).$$
On the other hand, $\bX \models R(a_0, \dots, a_{n-1})$, which contradicts the fact that $h$ is a homomorphism of $\bD$-separated structures.
\end{pf}

\separator

We now make a brief discussion of basic properties of the relation $\bea$.
First of all, we note that $\bea$ restricted to singletons is just a partial order.
Specifically, define
$$x \loe y \iff \sn x \bea \sn y.$$
By ($\i_2$) this relation is reflexive and antisymmetric and ($\i_3$) gives transitivity.
We shall come back to this observation later.
W shall call $\loe$ the \define{partial order associated to} $\bea$.
This is in fact the restriction of $\bea$ to pairs of one-element sets.
For some concrete classes of $\bea$-structures, it happens that the partial order $\loe$ is discrete, that is, $x \loe y$ holds only if $x = y$.

Another feature of $\bea$, in the presence of the constants $0$, $1$, is the ability to define the negation.
Namely, given a $\verazo$-separated structure $\bX = \pair X \bea$, we say that $b \in X$ is the \define{complement} (or \define{negation}) of an element $a \in X$ if both relations
$$\dn a b \bea \sn 0 \oraz \sn 1 \bea \dn a b$$
hold in $\bX$.
Using Proposition~\ref{PropTrzideset}, we see that $b$ is uniquely determined, because $\bX$ can be viewed as a family of sets $\Ef$, with $0$ being the empty set, $1$ the whole universe $V = \bigcup \Ef$, and then $\dn a b \bea \sn 0$ means that $a \cap b = \emptyset$ and $\sn 1 \bea \dn a b$ means that $a \cup b = V$.
On the other hand, it is an easy exercise to derive the uniqueness of the complement just from the axioms ($\i_0$)--($\i_3$) together with (\ref{EqC0}), (\ref{EqC1}).
The complement of $a$ will be denoted by $\neg a$.
In some cases, it will be natural to add the function $\neg$ to the language in order to get a duality.
This indeed happens with Stone duality.

We shall now try to explain the Pasch axiom.
Let us imagine that each finite set $x$ generates a ``lower-set" $\ell(x)$ and an ``upper-set" $u(x)$, possibly by adding elements outside of the universe $\abs \bX$.
Now the relation $a \bea b$ means that $u(a) \cap \ell(b) \nnempty$.
In fact, this is the case when $\bX = \pair \Ef \bea$ is like in Proposition~\ref{PropTrzideset}, where
$$\ell(x) = \setof{y}{y \subs \bigcup x} \oraz u(x) = \setof{y}{\bigcap x \subs y}.$$
Formally we may still have $\ell(a) \cap u(b) = \emptyset$ even though $\bigcap a \subs \bigcup b$, however as we said before, the lower- and upper-set may contain more elements, for instance, all sets that are finite unions of elements of $\Ef$. In that case, indeed $\ell(a) \cap u(b) \nnempty$ iff $a \bea b$.
Now let us imagine that both $\ell(x)$ and $u(x)$ mean the convex hull of the set $x$, with respect to some fixed (possibly different) convexity structures.
Condition ($\i_3$) becomes a geometric axiom saying that some convex
sets should intersect.
In case where $a_0, a_1, b_0, b_1$ are singletons, this is the classical Pasch axiom stated in a more abstract setting, namely that the line segments $[a_0, a_1]$ and $[b_0, b_1]$ should intersect whenever $a_1 \in [z, a_0]$ and $b_1 \in [z, b_0]$.

It turns out that a large subclass of $\vera$-separated structures indeed comes from convexity structures.
More details are explained in the next section, where we also review some classical dualities for 2-separated structures.

\section{Applications}\label{SecApplications}

In this section we study a duality involving convexity structures, as a special case of the duality involving $\bea$.
We also discuss some known dualities like Priestley duality for posets vs. distributive lattices and Hofmann-Mislove-Stralka duality for semilattices.

We already know that all 2-separated structures have their associated relation $\bea$ (formally a sequence of relations) which could replace the original language (see Theorem~\ref{Twthmeorgo}).
Given 2-separated classes $\Ef$ and $\Gee$, we shall say that $\Ef$ is \define{naturally dual} to $\Gee$ if for every finite $\bX \in \Ef$ it holds that $\bX^* \in \Gee$ and $\bX$ is reflexive with respect to the dual pair $\pair \vera \verazo$ or $\pair \wera \wera$ (depending on the existence of the constants 0, 1 in the language), and the same holds when interchanging the roles of $\Ef$ and $\Gee$.
In other words, $\Ef$ and $\Gee$ are naturally dual if the categories of finite structures of $\Ef$ and $\Gee$ are dually equivalent via the $\hom$-functors into $\vera$, $\wera$ or $\werazo$.
Note that we actually do not require that the classes $\Ef$ and $\Gee$ are of the form $\isp \bD$ and $\isp \bE$ for some dual pair $\pair \bD \bE$.
In any case, once we have suitable dual classes $\Ef$ and $\Gee$, we can use Theorems~\ref{tdwajedenxx} and \ref{tdwaipull} to conclude that all compact and all discrete structures in $\Ef$ and $\Gee$ are reflexive.

\subsection{Abstract convexities}\label{SubsEcconve}

A \define{convexity} on a set $X$ is, by definition, a family $\Gee$ of subsets of $X$ that is closed under arbitrary intersections and unions of chains.
So in particular $X \in \Gee$, being the intersection of the empty subfamily of $\Gee$.
Some authors assume that $\emptyset \in \Gee$, although we prefer to avoid it for the reasons explained below.
It is well known \cite{vandeVel} and easy to check that every convexity is determined by the \define{convex hull operator} defined by $\conv A$ to be the intersection of all convex sets containing $A$.
The relation ``$x \in \conv A$" can actually be defined in a first-order language, using infinitely many relations of the form $B(x,y_0, \dots, y_{n-1})$ meaning that $x \in \conv \{y_0, \dots, y_{n-1}\}$.

For the sake of generality, we need to work with two convexities on the same set.
Namely, a \define{bi-convexity space} is a structure of the form $\triple X \El \Yu$, where $\El$ and $\Yu$ are convexities on $X$.
We shall write $\conv_\El$ and $\conv_\Yu$, indicating which convexity we have in mind.
A bi-convexity space $\triple X \El \Yu$ is \define{normal} if
\begin{enumerate}[(N1)]
	\item For every $x \ne y$ in $X$ either $\conv_\El\sn x \cap \conv_\Yu \sn y = \emptyset$ or $\conv_\Yu\sn x \cap \conv_\El \sn y = \emptyset$.
	\item Given $A \in \El$, $B \in \Yu$, there exists $H \in \Yu$ such that $X \setminus H \in \El$ and $B \subs H$, $A \cap H = \emptyset$.
\end{enumerate}
It is easy to ``encode" normal bi-convexity structures in the language of $\bea$.
Namely, consider the following axiom:
\begin{enumerate}
	\item[($\i_4$)] $a \bea b \implies (\exists\; p) \;\; a \bea \sn p \Land \sn p \bea b$.
\end{enumerate}
Notice that the converse implication is a special case of ($\i_3$): put $a_0 = \emptyset$, $a_1 = a$, $b_0 = b$, $b_1 = \emptyset$.
The following fact is rather obvious.

\begin{prop}\label{Pnoita}
Let $\triple X \El \Yu$ be a normal bi-convexity space.
Given finite sets $a, b \subs X$, define
$$a \bea b \iff \conv_\Yu(a) \cap \conv_\El(b) \nnempty.$$
Then $\pair X \bea$ satisfies $(\i_0)$ -- $(\i_4)$.
\end{prop}

It turns out that the converse statement holds true, therefore we get an axiomatization of bi-convexity structures in terms of $\bea$.

\begin{prop}\label{Pnbotno}
Assume $\pair X \bea$ satisfies $(\i_0)$ -- $(\i_4)$.
Given a finite set $a \subs X$, define
$$\conv_\El(a) = \setof{p \in X}{\sn p \bea a} \oraz \conv_\Yu(a) = \setof{p \in X}{a \bea \sn p}.$$
Then $\triple X \El \Yu$ is a normal bi-convexity space satisfying the condition
\begin{equation}
a \bea b \iff \conv_\Yu(a) \cap \conv_\El(b) \nnempty.
\tag{$\ddagger$}\label{eqsztylet}
\end{equation}
\end{prop}

It is rather clear that the convexities $\El$ and $\Yu$ satisfying (\ref{eqsztylet}) are uniquely determined.

\begin{pf}
It is clear that $\conv_\El$ and $\conv_\Yu$ are monotone, that is, they induce convexities $\El$ and $\Yu$ on $X$.
We need to check that $\triple X \El \Yu$ is a normal bi-convexity space.
Condition (N1) is just the translation of ($\i_2$).
Fix disjoint sets $A \in \El$, $B \in \Yu$.
Given $s \in \fin B$, $t \in \fin A$, we have that $s \not \bea t$, therefore by the remark after the proof of Theorem~\ref{ThmPaschNormal} there exists a homomorphism $\map h X \vera$ such that $\img h A \subs \sn 0$ and $\img h B \subs 1$.
Then $h$ preserves the bi-convexity structures and hence $H = h^{-1}(1)$ is such that $H \subs B$, $H \cap A = \emptyset$ and $H \in \Yu$, $X \setminus H \in \El$.
This shows (N2).
\end{pf}

It is interesting to note that axiom $(\i_3)$ for bi-convexity spaces is equivalent to the following
\begin{align*}
q \in \conv_\Yu(a_0 \cup \sn p) & \Land r \in \conv_\El(b_1 \cup \sn p) \\
&\implies \conv_\Yu(a_0 \cup \sn r) \cap \conv_\El(\sn q \cup b_1) \nnempty,
\tag{$\P$}\label{eqPasch}
\end{align*}
which is indeed known as the Pasch axiom for abstract convexity structures (see \cite{Kub_sep} for more details).
Clearly, (\ref{eqPasch}) is a special case of $(\i_3)$.
On the other hand, it is easy to check that $(\i_3)$ follows from $(\i_0)$--$(\i_2)$,  $(\i_4)$ and (\ref{eqPasch}).

The two propositions above provide an axiomatization of normal bi-convexity structures in the language of $\vera$.
Note that $\vera$ itself can be regarded as a normal bi-convexity structure, because it satisfies $(\i_4)$ for the obvious reasons.
We are now ready to prove the main result of this subsection.

\begin{tw}\label{ThmConvexitiesn}
The class of normal bi-convexity structures is naturally dual to the class of normal bi-convexity structures with $0$ and $1$.
\end{tw}

\begin{pf}
We shall use Theorem~\ref{ThmDawidek}.
Fix a finite normal bi-convexity structure $\bX$.
We look at $\bX^*$ as the family of all subsets of $X$ whose characteristic functions are the homomorphisms from $\bX$ into $\vera$.
Fix $S, T \subs \bX^*$ such that $S \bea T$, that is, $\bigcap S \subs \bigcup T$.
Let $A = \bigcup_{t \in T}(X \setminus t)$.
Then $A$ is an $\El$-convex set disjoint from the $\Yu$-convex set $B = \bigcap S$.
Thus, by normality, there exists $H \in \Yu$ such that $X \setminus H \in \El$ and $B \subs H$, $A \cap H = \emptyset$.
The characteristic function of $H$ is a homomorphism into $\vera$, therefore $H \in \bX^*$, according to our agreement.
Finally, $S \bea \sn H$ and $ \sn H \bea T$, which shows that $\bX^*$ satisfies ($\i_4$).
By Theorem~\ref{ThmDawidek}, this completes the proof.
\end{pf}

It might be of some interest to describe structures dual to normal convexity structures, that is, normal bi-convexity structures $\triple X \El \Yu$ with $\El = \Yu$.
We do it below.

In view of Propositions~\ref{Pnoita} and \ref{Pnbotno}, the following axiom describes normal convexity spaces:
\begin{enumerate}
\item[($\i_5$)] $a \bea b \iff b \bea a$.
\end{enumerate}
A structure $\pair X \bea$ satisfying ($\i_5$) will be called \define{symmetric}.
As we know from Theorem~\ref{ThmConvexitiesn}, the dual to a normal convexity structure is necessarily a normal bi-convexity structure with $0$, $1$.
Notice that, given a normal convexity structure $\bX$, for every homomorphism $\map f \bX \vera$, the map $1-f$ is again a homomorphism.
It follows that the dual bi-convexity structure is closed under negation.
Note the following property:
\begin{equation}
a \bea b \iff \neg b \bea \neg a,
\tag{\ae}\label{EqNeggga}
\end{equation}
where $\neg s = \setof{ \neg p }{ p \in s }$.
We shall say that a bi-convexity space $\bX$ is \define{complemented} if $0$, $1$ are in the language of $\bX$ and $\bX$ is closed under negations.
If this is the case, then condition (\ref{EqNeggga}) is satisfied automatically.

\begin{tw}
The class of normal convexity structures is naturally dual to the class of complemented normal bi-convexity structures.
\end{tw}

\begin{pf}
Let $\bX$ be a finite normal convexity structure 
(i.e. a finite normal symmetric bi-convexity structure), %KP
represented in the language of $\vera$.
By Theorem~\ref{ThmConvexitiesn}, the dual $\bX^*$ is a normal bi-convexity structure 
with $0$ and $1$.  %KP
Clearly, it is complemented.
Conversely, if $\bX$ is a complemented normal bi-convexity structure then, again by Theorem~\ref{ThmConvexitiesn}, the dual $\bX^*$ is a normal bi-convexity structure.
Since $\bX$ is complemented, $\bX^*$ is symmetric, that is, both convexities are the same.
\end{pf}

Note that the class of $2$-separated convexity spaces (or, more generally, bi-convexity spaces) can be described as $\isp \bG$, where $\bG$ is the unique normal (bi-)convexity structure whose universe is $\dn 01$.
Recall that the language of $\bG$ should consist of all relations coding ``$x \in \conv \{y_1, \dots, y_n\}$".
Unfortunately, the class $\isp \bG$ contains convexity structures that are not normal.
In fact, it is easy to find a 5-element subset in the plane that is not normal as a convexity space.

We do not know whether there exists a structure $\bE$ such that $\pair \bG \bE$ is a dual pair.

\subsection{The Hofmann-Mislove-Stralka duality}

Fix a finite semilattice $\bX = \pair X \meet$, that is, $\meet$ is a symmetric associative operation on $X$ such that $x \meet x = x$ holds for every $x \in X$.
We use the adjective ``meet" only to indicate that $\meet$ plays the role of the infimum with respect to the associated partial order defined by
$$x \loe y \iffdef x \meet y = x.$$
Given finite sets $s, t \subs X$, define
$$s \bea t \iffdef \Meet s \loe t,$$
where $\loe$ is the associated partial order.
It is very easy to check that $\bea$ satisfies ($\i_0$)--($\i_3$) and therefore $\bX$ becomes a $\vera$-separated structure.
Actually, $\bea$ is the same as the relation defined by condition (\ref{Eqwationae}) before Theorem~\ref{Twthmeorgo}.

Given $f \in \hom(\bX,\vera)$, observe that the set $F = f^{-1}(1)$ is meet-closed, therefore it is of the form $[p,\rightarrow)$, where $p = \Meet F$.

For a moment, let us identify each $f \in \hom(\bX,\vera)$ with the set $f^{-1}(1)$.
Then, given finite sets $S, T \subs \hom(\bX, \vera)$ we have that $S \bea T$ if and only if $\bigcap S \subs T$, just because of the remark above.
It follows that the dual to a semilattice (with respect to the $\pair \vera \verazo$ duality) is again a semilattice.

By Theorem~\ref{TwCorWeronisia}, we obtain the Hofmann-Mislove-Stralka duality:

\begin{tw}[Hofmann, Mislove, Stralka \cite{HMS}]
The class of semilattices is naturally dual to the class of semilattices with 0 and 1.
\end{tw}

We also get another variant:

\begin{wn}
The class of semilattices with 0 is naturally dual to itself.
\end{wn}

\subsection{Priestley and Stone dualities}

We already know that every $\vera$-separated structure has a natural partial order, namely, $\bea$ restricted to pairs of singletons.
In some cases, this partial order trivializes, that is, it is just equality.
This happens, for example, in symmetric structures.
Let us now consider the other extreme: The case where the partial order induces the relation $\bea$.
The appropriate axiom is:
\begin{equation}
s \bea t \iff (\exists \; p \in s)(\exists \; q \in t) \;\; p \loe q.
\tag{$\leq$}\label{Eqposetblrjjr}
\end{equation}
The class of $\vera$-separated structures satisfying (\ref{Eqposetblrjjr}) coincides with 
$\isp{\bO}$, 
where $\bO = \pair 2 \loe$ is the 2-element linearly ordered set in which $0 < 1$.

Recall that a \define{bounded distributive lattice} is a structure of the form $\seq {L, \meet, \join, 0, 1}$ isomorphic to a family of sets $\seq{\El, \cap, \cup, \emptyset, X}$, where $\bigcup \El = X$.

It is straight to see that the dual to a $\bO$-separated structure is closed under meet and join.
It also carries the constants $0$, $1$, therefore it is a bounded distributive lattice.

Conversely, given a bounded distributive lattice $\bL$, its dual structure (with respect to the $\pair \vera \verazo$ duality) is a partially ordered set.
This follows from the following folklore fact:

\begin{prop}
Let $S$, $T$ be finite families of prime filters in a distributive lattice such that $\bigcap S \subs \bigcup T$.
Then there exist $p \in S$ and $q \in T$ such that $p \subs q$.
\end{prop}

Summarizing, we obtain Priestley duality as a special case of Theorem~\ref{TwCorWeronisia}:

\begin{wn}[Priestley]
Partially ordered sets are naturally dual to bounded distributive lattices.
\end{wn}

Stone duality can actually be regarded as a special case of Priestley duality.
Namely, trivial partial orders correspond to complemented distributive lattices, i.e., Boolean algebras.
In other words:

\begin{wn}[Stone]
Sets (i.e. structures whose language consists of equality only) are naturally dual to Boolean algebras.
\end{wn}

Note that Priestley duality provides a simple proof of the classical fact saying that a bounded distributive lattice is a Boolean algebra if and only if each of its prime filters is maximal.
Of course, using Theorems~\ref{tdwajedenxx} and \ref{tdwaipull}, we obtain the full Stone duality, saying that compact 0-dimensional spaces are dually equivalent to Boolean algebras.

Since the duality is symmetric, one may ask for compact Boolean algebras and their duals, which should be just sets (with no structure). Obviously, we get the class of all Cantor cubes. In fact, this is again easily proved using duality:

\begin{prop}
Let $K$ be a compact Boolean algebra which is $\dn 01$-separated. Then $K$ is topologically isomorphic to $2^\kappa$ for some cardinal $\kappa$.
\end{prop}

\begin{pf}
Let $A = K^*$.
Recall that $A$ carries no structure, therefore every function $\map fA2$ is a homomorphism.
By Theorem~\ref{tdwajedenxx}, $K = 2^A$.
\end{pf}

We finish with an easy example showing that semi-dual pairs need not be dual.

\begin{ex}
As before, let $\bL$ be the 2-element bounded distributive lattice and let $\bD = \pair 2=$.
We claim that $\pair \bD\bL$ is a semi-dual pair.
Indeed, if $\X\in\isp\bD$ then $\X^* = \hom(\X,\bD)$ has a structure of a distributive lattice (in fact: it is a Boolean algebra).
If $\X\in\isp\bD$ is finite then $\X = \X^{**}$, because of the Priestley duality.
On the other hand, if $\Y \in \isp\bL$ is a finite distributive lattice which is not a Boolean algebra then $\Y\ne \Y^{**}$, because $P^{*}$ is a Boolean algebra whenever $P$ is a $\bD$-separated structure (i.e. $P$ is any set with equality).
It follows that $\pair \bD\bL$ is not a dual pair.
\end{ex}

\subsection{Reflexivity of topological posets}

It turns out that many topological $2$-separated posets are reflexive.
Below we give a necessary and sufficient condition for reflexivity.
Given a poset $P$ and $A,B\subs P$ we write $A \loe B$ if $x \loe y$ holds for every $x\in A$, $y\in B$.
We write $\fs(P)$ for the set of all \define{final segments} in $P$, that is, all sets $F\subs P$ whose characteristic functions are order preserving.
Clearly, $\fs(P)$ can be identified with the dual of $P$.

\begin{lm}\label{sdgsdgwe}
Let $\poset=\triple P\loe\Tau$ be a $2$-separated topological poset. Let $F$ be a closed filter in $\poset^*$ and define $B=\bigcap F$. Then
$$F=\setof{x\in\fs(\poset)}{B\subs x}.$$
In particular, $B\nnempty$ (recall that, by definition, $\emptyset\notin F$).
\end{lm}

\begin{pf}
Fix $x_0\in\fs(\poset)$ such that $B\subs x_0$. It suffices to show that $x_0\in\cl B$. For this aim, fix a basic neighborhood of the form $s^+\cap t^-$ of $x_0$, where $s,t\in\fin P$. Since $t\cap B=\emptyset$ and $F$ is closed under finite intersections, there is $y\in F$ such that $y\cap t=\emptyset$. Let $x=x_0\cup y$. Then $x\sups x_0$, therefore $x\in F$. Clearly $x\in s^+\cap t^-$.
\end{pf}

\begin{tw}
Let $\poset=\triple P\loe\Tau$ be a $2$-separated topological poset in which all clopen final and initial segments generate the topology $\Tau$. Then $\poset$ is reflexive if and only if for every closed final segment $B\subs P$ and for every closed initial segment $C\subs P$ with $C<B$, there exists a clopen final segment $u$ such that $B\subs u$ and $C\cap u=\emptyset$. 
\end{tw}

\begin{pf}
We start with the ``if" part. Fix a poset $\poset$ satisfying the above condition.
Let $H$ be a clopen prime filter and let $B=\bigcap H$. By Lemma \ref{sdgsdgwe}, $u\in H$ iff $B\subs u$ and $u\in\fs(\poset)$. Let $I=\fs(\poset)\setminus H$. Then $I$ is a clopen (prime) ideal. Let $C=P\setminus\bigcup I$. We claim that
$$I=\setof{u\in\fs(\poset)}{u\cap C=\emptyset}.$$
In fact this statement follows from Lemma \ref{sdgsdgwe}, considering the dual lattice of initial segments.
Now we claim that $C\loe B$. Indeed, suppose $x\in C$, $y\in B$ are such that $x\not\loe y$ and find $u\in\fs(\poset)$ with $x\in u$ and $y\notin u$. Then $u\notin I$ and $u\notin H$, a contradiction. Since $I\cap H=\emptyset$, we see that $B$ and $C$ cannot be separated by a clopen final segment. By our assumption, this means that $B\cap C\nnempty$. Finally, if $q\in B\cap C$ then $B=[q,\rightarrow)$ and hence $H=q^+$.

Now suppose that $C,B$ witness the failure of the condition in the lemma. Define $H=\setof{u\in\fs(\poset)}{B\subs u}$ and $I=\setof{u\in\fs(\poset)}{u\cap C=\emptyset}$. Then $H$ is a closed filter, $I$ is a closed ideal and $I\cap H=\emptyset$ since $u\in I\cap H$ would separate $C$ and $B$. Moreover $I\cup H=\fs(\poset)$ because if $u\notin I$ then there is $p\in C\cap u$ and therefore $B\subs [p,\rightarrow)\subs u$, that is $u\in H$. Thus $H$ is a clopen prime filter.
We claim that $H\ne q^+$ for any $q\in P$. Indeed, suppose $H=q^+$; then $q\notin C$ so there exists a neighborhood $v$ of $q$ which is disjoint from $C$. We may assume that $v=u\setminus w$, where $u,w\in\fs(\poset)$. Then $q\in u$, so $B\subs u$. Thus $C\cap u\nnempty$, because $u$ does not separate $B$ from $C$. Now if $p\in C\cap u$ then $p<q$ and hence $p\in u\setminus w$, because $P\setminus w$ is an initial segment. Thus $v\cap C\nnempty$, a contradiction. This completes the proof of the ``only if" part.
\end{pf}

Below is an example of an irreflexive poset.

\begin{ex} Define $S=\sn0\cup\setof{1/n}{\ntr}$ with the natural linear order and let
$X=(\omega_1+1)\times S$ be endowed with the product ordering $\loe$ and with the product topology. Clearly, $X$ is $2$-separated. Define
$\poset=\triple{P}{\loe_P}\Tau$, where $P=X\setminus\sn{\pair{\omega_1}{0}}$, $\loe_P=\loe\rest (P\times P)$ and $\Tau$ is the subspace topology. It is clear that $\poset$ is $2$-separated. Let $B=\setof{\pair{\omega_1}{1/n}}{\ntr}$ and define
$$H=\setof{u\in\fs(\poset)}{B\subs u}.$$
Clearly, $H$ is a closed filter in $\poset^*$. 
Observe that $B=\bigcap H$ and $H\ne p^+$ for any $p\in P$. We shall show that $H$ is a clopen prime filter.
Let $C=\omega_1\times\sn0$ and define
$$I=\setof{u\in\fs(\poset)}{u\cap C=\emptyset}.$$
Clearly, $I$ is a closed ideal in $\poset^*$. Now it suffices to show that $I\cup H=\fs(\poset)$ and that $I\cap H=\emptyset$.
Fix $u\in \fs(\poset)\setminus I$. Then there is $p\in C\cap u$ and hence $B\subs [p,\rightarrow)\subs u$, which shows that $u\in H$. Hence $I\cup H=\fs(\poset)$.
Fix $u\in H$. Then $u$ is open and $B\subs u$. Thus for each $\ntr$ there is $\delta(n)<\omega_1$ with $[\delta(n),\omega_1]\times\sn{1/n}\subs u$. Let $\delta=\sup_{\ntr}\delta(n)$.
Then $\delta<\omega_1$ and we have $\pair\delta0\in\cl u=u$. It follows that $u\cap C\nnempty$, that is $u\notin I$. This shows that $I\cap H=\emptyset$.
\end{ex}

\subsection{Betweenness structures}

There exist quite natural $2$-element structures for which there is no duality.
Below we describe one.
Let $\beto$ be the following ternary relation on $2$:
$$\beto(x,y,z) \equiv (x = z = 1 \implies y = 1).$$
Let $\Szero = \pair 2\beto$.
An $\Szero$-separated structure is called a \define{betweenness structure}, sometimes also called $S_0$-betweenness.
A subset $G$ of a betweenness structure $\pair XB$ is \define{convex} if
$$(\forall\; a,b \in G) (\forall\; x\in X)\; B(a,x,b) \implies x\in G.$$
We also define $[a,b] = \setof{x\in X}{B(a,x,b)}$, sometimes called the \define{interval} joining $a$ and $b$.
It is clear that convex sets form a (usually non-normal) convexity structure in the sense of Subsection~\ref{SubsEcconve}, however not all convexities are induced by a betweenness relation.

\begin{lm}
A structure $\pair XB$ is $\Szero$-separated if and only if $B$ satisfies the following axioms.
\begin{enumerate}
	\item[(1)] $B(x,x,y)$ and $B(x,y,y)$.
	\item[(2)] $B(u,x,v) \Land B(u,y,v) \Land B(x,z,y) \implies B(u,z,v)$.
	\item[(3)] $B(x,y,x) \Land B(y,x,y) \implies x = y$.
\end{enumerate}
If $\pair XB$ is a betweenness structure then for every $a,b\in X$ the set $[a,b]$ is convex.
\end{lm}

\begin{pf}
It is rather clear that every structure in $\isp\Szero$ satisfies the above axioms.
Fix $\pair XB$ such that $B$ satisfies (1) -- (3).
Notice that (1) and (2) imply the second statement, i.e. that the set $[a,b]$ is convex for every $a,b\in X$.
Now if $a\ne b$ then by (3) either $a\notin [b,b]$ or $b\notin [a,a]$.
Further, the characteristic function of every convex set is a homomorphism.
It follows that $\pair XB$ is $\Szero$-separated.
\end{pf}

\begin{tw}
There exist no structure $\bE = \pair 2\R$ for which $\pair \Szero \bE$ is a dual pair.
\end{tw}

\begin{pf}
Suppose $\bE$ is such a structure.
Let $\X = \pair \omega B$, where
$$B(k,\ell,m) \iff (k = m \implies \ell = k = m).$$
This is in some sense the ``minimal" betweenness structure on $\omega$.
Namely, non-trivial convex sets are the singletons of $\omega$.
Let $\map\Eva\X{\X^{**}}$ be the evaluation map.
By Theorem~\ref{tdwaipull}, we know that $\img\Eva\X = \X^{**}$.
On the other hand, $\X^{**}$ is pointwise closed in the space of continuous functions $C(\X^*,2)$.

Every homomorphism from $\X$ into $\Szero$ is the characteristic function of a convex subset of $\omega$.
Thus, $\X^*$ is naturally homeomorphic to $K = \dn\emptyset\omega \cup \dpower\omega1$.
After this identification, we have that
\begin{equation}
\Eva(n)^{-1} = \dn \omega n
\tag{$\dagger$}\label{eqwgowrgpj}
\end{equation}
for every $\ntr$.
The set $\setof{\Eva(n)}{\ntr}$ should be pointwise closed.
On the other hand, (\ref{eqwgowrgpj}) shows that the characteristic function of $\sn \omega$ is the limit of the sequence $\sett{\Eva(n)}{\ntr}$.
This is a contradiction, because $\omega$ is an isolated point of $K$ and therefore its characteristic function is continuous.
\end{pf}

\begin{ex}
Let $B$ be the natural betweenness relation on $\bD = \dn01$. That is, $B(x,y,z)$ holds iff $y \in \dn xz$.
Let us describe the class $\isp\bD$.

Fix a set $X$ and a ternary relation $B$ on $X$. We say that $B$ is a \define{betweenness} relation if the following conditions are satisfied for every $x,y,z,a,b \in X$.
\begin{enumerate}
	\item[(1)] $B(x,x,y)$ and $B(x,y,y)$.
	\item[(2)] If $B(a,x,b)$, $B(a,y,b)$ and $B(x,z,y)$ then $B(a,z,b)$.
\end{enumerate}
It is straight to check that every $\pair XB \in \isp \bD$ satisfies the above conditions.
On the other hand, if $B$ is a betweenness relation on a set $X$ then $\pair XB \in \isp \bD$ if and only if for every $x,y,z\in X$ such that $\neg B(x,z,y)$ there exists a halfspace $H\subs X$ such that $z\notin H$ and $x,y \in H$.
Recall that a set $H$ is a halfspace if its characteristic function is a homomorphism into $\bD$.
In other words, $H$ is a halfspace iff whenever $B(x,z,y)$ holds and $x,y\in H$ then $z\in H$ and whenever $B(x,z,y)$ holds and $x,y\notin H$ then $z\notin H$.

\begin{pyt}
Does there exist $\bE$ such that $\pair \bD\bE$ is a dual pair?
\end{pyt}
\end{ex}

\section{Notational remarks}\label{SecNotationalremarks}  %KP

In this section we assume that the reader knows the book \cite{CD} (or another texts on algebraic dualities) and our intention here is to explain the key similarities and differences between the approach and terminology in this paper and the above-cited book. Readers who are not familiar with the theory of algebraic dualities can safely stop reading here.

Perhaps the key difference with the approach in \cite{CD} is our definition of semi-dual pair which just describes the conditions for duality between finite structures.
%We don't try to formulate (algebraic or another) conditions on how should $\bE$ relate to $\bD$ to obtain a semi-dual pair. We feel that at this level of generality speaking directly what we want (see the definition of semi-dual pair) is the best way.
As it was signalized earlier, here we use the term \textit{natural duality} as something similar to the full duality for the class of finite structures. 

Because we know that the duality for finite structures implies the duality for both compact and discrete structures (see Section~\ref{SecNatDual}), we decided to use the property of finite structures in the definitions.

\newcommand{\Mf}{\ensuremath{\mathrm{M}\kern-0.75em\raise-.9ex\hbox{$\sim$}}}
\newcommand{\Mp}{\ensuremath{\underline{\mathrm{M}}}}

Observe that if \Mf\ yields a natural duality (in the sense of \cite[Chapter~2]{CD}) on $\isp{\Mp}$ %and \Mf\ does not contain partial operations,
then the $\pair \Mp\Mf$ is a semi-dual pair.
Also note that since we allow more general structures in the definition of semi-dual pairs, the inverse implication does not need to hold (or at least is hard to express in the language of algebraic dualities). % Directly from 
% In other words, if 

Observe that Theorem \ref{tdwaipull} addresses a problem similar to those presented in Duality Compactness Theorem (see \cite[2.2.11]{CD}) and Third Strong Duality Theorem (see \cite[3.2.11]{CD}).

The family of relations $\bea$ (see subsection \ref{SubSecUlt}) not only essentially ``refines'' the Brute Force construction (see \cite[section 2.1]{CD}) for the two-element schizophrenic object, but also may be applied to full duality. 

Note that our results may be applied for more general classes than $\isp{\bE}$; see Section \ref{Sec2sep} and the beginning of Section \ref{SecApplications}. 
Finally, some of the results of the last section cannot be obtained by classical methods of the theory of algebraic dualities.

%Our proofs of classical dualities 

% Probably the biggest difference is in our arguments: we use more pure set theory and topology, less algebra.

\printindex


\begin{thebibliography}{99}

\bibitem{ABKR} {\sc Abraham, U.; Bonnet, R.; Kubi\'s, W.; Rubin, M.}, {\em On poset Boolean algebras\/}, Order 20 (2003), no. 3, 265--290. 

\bibitem{Banasch} {\sc Banaschewski, B.},
{\em Remarks on dual adjointness\/}, Nordwestdeutsches Kategorienseminar (Tagung, Bremen, 1976), pp. 3–10. Math.-Arbeitspapiere, No. 7, Teil A: Math. Forschungspapiere, Univ. Bremen, Bremen, 1976. 

%\bibitem{Birkhoff} {\sc Birkhoff, G.}, {\em Lattice theory\/}. Third edition. American Mathematical Society Colloquium Publications, Vol. XXV American Mathematical Society, Providence, R.I. {\bf 1967}.

\bibitem{CD} {\sc Clark, D.M.; Davey, B.A.}, {\em Natural Dualities for the Working Algebraist\/}. Cambridge Studies in Advanced Mathematics, 57. Cambridge University Press, Cambridge, 1998.

% \bibitem{CDJP} {\sc Clark, D.M.; Davey, B.A.; Jackson, M.G.; Pitkethly, J.G.}
% {\em The axiomatizability of topological prevarieties}, 
% Advances in Mathematics 218 (2008), no. 5, 1604--1653.

\bibitem{Davey2007} {\sc Davey, B.A.}, {\em Natural dualities for structures\/}, in: 
{\em Acta Universitatis Matthiae Belii Series Mathematics}, pp. 3--28, Banska Byst\v rica {2007}. %KP

\bibitem{DHP} {\sc Davey, B.A.; Haviar, M.; Priestley, H.A.}, {\em Natural dualities in partnership\/}, Appl. Categor. Struct. 20 (2012) 583--602.

\bibitem{DP} {\sc Davey, B.A.; Priestley, H.A.}, {\em Introduction to Lattices and Order\/}, Cambridge University Press {2002}.

%\bibitem{Engelking} {\sc Engelking, R.}, {\em General Topology\/}. Translated from the Polish by the author. Second edition. Sigma Series in Pure Mathematics, 6. Heldermann Verlag, Berlin, {\bf 1989}. %KP

%\bibitem{Gratzer} {\sc Gr\"atzer, G.}, {\em General Lattice Theory\/}. Second edition. Birkh\"auser Verlag, Basel, {\bf 1998}.

\bibitem{Hof} {\sc Hofmann, D.}, {\em A generalization of the duality compactness theorem}, J. Pure Appl. Algebra 171 (2002), 205--217. %KP

\bibitem{HMS} {\sc Hofmann, K.H.; Mislove, M.; Stralka, A.}, {\em The Pontryagin Duality of Compact ${0}$-dimensional Semilattices and its Applications\/}, 
Lecture Notes in Mathematics, Vol. 396, Springer-Verlag, Berlin-New York, 1974.

%\bibitem{Is} {\sc Isbell, J.R.}, {\em Median algebra\/}, Trans. Amer. Math. Soc. 260 (1980), no. 2, 319--362.

\bibitem{Isbell} {\sc Isbell, J.R.}, {\em General functorial semantics, I\/}, Amer. J. Math. 94 (1972) 535--596.

\bibitem{Jo} {\sc Johansen, S.M.}, {\em Natural dualities for three classes of relational structures}, Algebra Univers. 63 (2010), no. 2--3, 149--170. %KP

\bibitem{Kub_sep} {\sc Kubi\'s, W.}, {\em Separation properties of convexity spaces\/}, J. Geom. 74 (2002) 110--119.

\bibitem{vandeVel} {\sc van de Vel, M.}, {\em Theory of Convex Structures\/}, North-Holland Mathematical Library 50, Amsterdam 1993.

\end{thebibliography}
\end{document}